\documentclass[11pt]{article}
\usepackage{amsfonts}
\usepackage{latexsym}
\usepackage{color, amsmath,amssymb, amsfonts, amstext,amsthm, latexsym, dsfont}
\usepackage{epic,eepic}
\usepackage{graphicx}
\usepackage{epstopdf}
\usepackage{longtable}
\usepackage[title]{appendix}
\usepackage{indentfirst}
\usepackage{graphicx}
\usepackage{float}
\usepackage{graphicx}
\usepackage[utf8]{inputenc}
\usepackage{subfigure}
\usepackage{url}
\usepackage{booktabs}
\usepackage[colorlinks=false, linktocpage=true]{hyperref}
\usepackage{caption}
\allowdisplaybreaks[4]

\usepackage{tikz}
\usepackage{etoolbox}
\usepackage{enumitem}

\definecolor{awesome}{rgb}{1.0, 0.13, 0.32}
\definecolor{bronze}{rgb}{0.8, 0.5, 0.2}
\definecolor{burntorange}{rgb}{0.8, 0.33, 0.0}
\definecolor{blue(ncs)}{rgb}{0.0, 0.53, 0.74}
\theoremstyle{plain}

\theoremstyle{remark}

\theoremstyle{definition}

\newcommand{\B}{\bf \color{blue(ncs)} }

\numberwithin{equation}{section} 

\begin{document}

\title{\bf Detecting the Most Probable High Dimensional Transition Pathway Based on Optimal Control Theory}
\author{\bf\normalsize{
Jianyu Chen$^{1,}$\footnote{Email: \texttt{jianyuchen@hust.edu.cn}}, Ting Gao$^{1,}$\footnote{Email: \texttt{tgao0716@hust.edu.cn}}, Yang Li$^{2,}$\footnote{Email: \texttt{liyangbx5433@163.com}}, Jinqiao Duan$^{3,}$\footnote{Email: \texttt{duan@iit.edu}}}\\[10pt]
\footnotesize{$^1$Center for Mathematical Sciences, Huazhong University of Science and Technology,} \\
\footnotesize{$^2$School of Automation, Nanjing University of Science and Technology, Nanjing 210094, People’s Republic of China} \\
\footnotesize{$^3$Departments of Applied Mathematics \& Physics, Illinois Institute of Technology,}\\
\footnotesize{Chicago, IL 60616, USA.}\\
}

\date{}
\maketitle
\vspace{-0.3in}
\begin{abstract}
Many natural systems exhibit phase transition where external environmental conditions spark a shift to a new and sometimes quite different state. Therefore, detecting the behavior of a stochastic dynamic system such as the most probable transition pathway, has made sense.

We consider stochastic dynamic systems driven by Brownian motion. Based on variational principle and Onsager-Machlup action functional theory, the variational problem is transformed into deterministic optimal control problem. Different from traditional variational method, optimal control theory can handle high dimensional problems well. This paper describes the following three systems, double well system driven by additive noise, Maire-Stein system driven by multiplicative noise, and Nutrient–Phytoplankton–Nooplanktonr system.

For numerical computation, based on Pontryagin's Maximum Principle, we adopt method of successive approximations to compute the optimal solution pathway, which is the most probable transition pathway in the sense of Onsager-Machlup. 
     
    \textbf{Keywords and Phrases:} Most probable transition pathway, optimal control, Pontryagin's Maximum Principle, method of successive approximations
\end{abstract}
\section{ Introduction}

Optimal control theory is a nontrivial extension of the classical theory of calculus of variations in two main directions: to dynamical and non-smooth settings. This builds on important contributions of Weierstrass and others and led in two inter-related directions: the Pontryagin’s maximum principle and the Hamilton-Jacobi-Bellman theory. An interesting historical account of the developments can be found in the book \cite{Liberzon2013Calculus}.

Random fluctuations in dynamical systems may lead to unexpected rare events.  
Transitions of a multistable system will reflect dramatic changes in the real world. Stochastic differential equations are widely used in various areas, including physics \cite{1997Infinite}, biology \cite{2012Experimental}, engineering and finance \cite{2003Mathematical}. They take environmental random fluctuations into consideration, which makes themselves important models in simulating real phenomena and predicting rare events.

In the presence of the noise, the dynamical behavior described by the stochastic differential equation is totally different compared to the deterministic differential equations. \cite{2023The} Transitions between equilibrium states of the vector field are prohibited for deterministic differential equations. However, no matter how small the noise is, there is transition between equilibrium states of the vector field for stochastic differential equations.

The major part of literature studying stochastic differential equations is devoted to the study of Gaussian dynamics, that is, the stochastic differential equations driven by Brownian motion. This has been applied in various fields. Biswas et.al focuses on characterising numerically the attractor volume in the state space of dynamical systems excited by additive white noise \cite{2021Characterising}. Wang investigates the regime shift behaviors between the prosperous state and the extinction state and stochastic resonance phenomenon for a meta-population system subjected to time delay and correlated Gaussian colored noises \cite{2021Impact}. It can be seen that dynamic systems driven by Gaussian noise are very common in nature, and a large number of literatures show that this is of practical significance.

Optimization problem is a derivative of dynamical system. One of the principal approaches in solving optimization problems is to derive a set of necessary conditions that must be satisfied by any optimal solution. Optimal control problems may be regarded as optimization problems in infinite-dimensional spaces, thus they are substantially difficult to solve. Although sufficient and necessary conditions for first-order and second-order optimization are proposed in \cite{Liberzon2013Calculus}, it still leaves great challenges for numerical computing. The maximum principle, formulated and derived by Pontryagin and his group in the 1950s, is truly a milestone of optimal control theory. It states that any optimal control along with the optimal state trajectory must solve the so-called (extended) Hamiltonian system, \cite{Liberzon2013Calculus} which is a two-point boundary value problem (and can also be called a forward-backward differential 
equation, to be able to compare with the stochastic case), plus a maximum condition of a function called the Hamiltonian. The mathematical significance of the maximum principle lies in that maximizing the Hamiltonian is much easier than the original control problem that is infinite-dimensional. This leads to closed-form solutions for certain classes of optimal control problems, including the linear quadratic case. The maximum principle has demonstrated its application in various disciplines. Bartholomew-Biggs \cite{0The} has done the optimisation of spacecraft orbital manoeuvres with Pontryagin's maximum principle. 

{\B Background of Onsager-Machlup functional and problem in computing the optimal transition pathway.}
The Onsager-Machlup functional is a powerful tool in tackling with transition behavior. Roughly speaking, Onsager-Machlup functional principle measures the maximum probability of rare events, which is a minimizer of the action of pathways connecting two ends. The minimizer is regarded as the most probable transition pathways in a wide range of literature. Hu et al \cite{hu2021transition}
have already applied Onsager-Machulp functional theory to compute the most probable transition pathway for a noised-induced stochastic system.

In this paper, we aim to compute the most probable transition pathways for stochastic systems with Pontryagin's maximum principle. The commonly used variational method can handle this problem. Hu et al \cite{hu2021transition} uses the variational method to calculate the optimal orbit by deducing Euler-Lagrange equation and combining with the neural network shooting method. It was proposed to compute the most probable transition pathway for stochastic differential equations with metastable states. The basic idea based on variational principle is to find the corresponding Euler-Lagrange equation of the system and solve it. However, when it comes to high dimensional systems, the problem becomes a little tricky. We have been looking for the alternatives.

{\B Background of Pontryagin's Maximum Principle}. Optimal control theory is a nontrivial extension of the classical theory of calculus of variations. This builds on important contributions of Weierstrass and others and led in two inter-related directions: the Pontryagin’s maximum principle and the Hamilton-Jacobi-Bellman theory. An interesting historical account of the developments can be found in \cite{Liberzon2013Calculus}.

The maximum principle was developed by the Pontryagin school in the Soviet Union in the late 1950s. \cite{Liberzon2013Calculus} It was presented to the wider research community at the first IFAC World Congress in Moscow in 1960 and in the celebrated book. It is worth reflecting that the developments we have covered so far in this book \cite{Liberzon2013Calculus}—starting from the Euler-Lagrange equation, continuing to the Hamiltonian formulation, and culminating in the maximum principle—span more than 200 years. The progress made during this time period is quite remarkable, yet the origins of the maximum principle are clearly traceable to the early work in calculus of variations.

The constrained minimization problem associated with computation of the most probable transition pathway can be naturally regarded as an optimal control problem. There is a huge literature of methodologies in the field of optimal control theory, involving time or space-marching \cite{1995The}, neural network \cite{2005Nearly}, stochastic process, Hopf formula \cite{2004An} and et al. For the stochastic systems driven by Brownian motion, the most probable transition pathway is directly the solution of a general optimal control problem according to Onsager-Machlup functional. 

{\B Introduce the machine learning method and its advantages.} In recent years, the machine learning method has greatly influenced areas involving computing. They can be applied to numerical solve partial differential equations and get rid of curse of dimensionality. They can be used to accelerating the process to obtain the minimizer of an optimal control problem. Applications to auto-driving \cite{2010Using}, auto-coding and even auto-proving. We take the advantages of machine learning method to our problem in finding the most probable transition pathway.

{\B Brief introduction of our methods and results.} We formulate the problem in finding the most probable transition pathway as a constrained minimization problem. The training procedure of a neural network can be naturally used in finding the constrained minimizer when a proper loss function is constructed.

{\B Arrangement of the article.}
This paper is arranged as follows. In section \ref{sec2}, we first briefly introduce the Onsager-Machlup theory and associated minimization problem in finding the most probable transition pathway in subsection \ref{OM}. Then we convert such problems from variational into optimal control view \ref{optimal_control_problem}. Last we state the general form of Pontryagin's Maximum Principle in subsection \ref{PMP}. Particularly attention is given to construct the algorithm for numerically solving our problem in section \ref{Algorithm}. Numerical computation will be given in section \ref{application}.


\section{From Calculus of Variations to Optimal Control}\label{sec2}

In the view of calculus of variations, we have treated so far problems with minimizing a cost functional over a given family trajectories. While then, optimal control theory studies similar problems but from a more dynamic viewpoint.

\subsection{Calculus of Variations Based Method for Detecting the Most Probable Transition Pathway} \label{OM}

A classical result \cite{Liberzon2013Calculus} due to Euler and Lagrange gives a necessary condition for optimality that allows us to solve such problem. We briefly introduce the Onsager-Machlup theorem and the most probable transition pathway.

We consider a generalized stochastic differential equation with multiplicative Gaussian noise in $\mathbb{R}^d$:
\begin{equation}\label{generalSDE}
\begin{split}
dX(t)=\tilde{b}(X(t))dt+\sigma(X(t))dB(t), t\in[0,t_{f}],
\end{split}
\end{equation}
with initial condition $X(0)=x_0 \in \mathbb{R}^d$, where $\tilde{b}:\mathbb{R}^d\rightarrow\mathbb{R}^d$ is a regular function, $\sigma:\mathbb{R}^d\rightarrow\mathbb{R}^{d\times k}$ is a $d\times k$ matrix-valued function, and $B$ is a Brownian motion in $\mathbb{R}^k$. The well-posedness of stochastic differential equation \eqref{generalSDE} has been widely investigated \cite{karatzas2012brownian}. We are interested in the transition phenomena between two meta-stable states (taking as the stable equilibrium states of the corresponding deterministic system) if there exist. 

To investigate transition phenomena, one should estimate the probability of the solution pathways on a small tube. The Onsager-Machlup theory of stochastic dynamical system 
 \eqref{generalSDE} gives an approximation of the probability
\begin{equation}\label{probability}
\begin{split}
\mathbb{P}(\{\|X-z\|_{t_{f}} \leqslant \delta\}) \propto C(\delta,t_{f}) \exp \left\{-S(z,\dot{z})\right\},
\end{split}
\end{equation}
where $\delta$ is positive and sufficiently small, $C(\delta,t_{f})=\mathbb{P}(\{\|B\|_{t_{f}} \leqslant \delta\})$, the probability of Bownian motion in the $\delta$-tube, $\|\cdot\|_{t_{f}}$ is the uniform norm of the space of all continuous functions in the time interval $[t_{0},t_{f}]$, $z$ is function in this space and the Onsager-Machlup action functional is

\begin{equation}\label{action}
\begin{split}
S(z,\dot{z})=\frac{1}{2}\int_{0}^{t_{f}}[(\dot{z}-b(z))V(z)(\dot{z}-b(z))+\operatorname{div}b(z)-\frac{1}{6}R(z)]dt,
\end{split}
\end{equation}
where $V(z)=(\sigma(z)\sigma^{*}(z))^{-1}$, $R(z)$ is the scalar curvature with respect to the Riemannian metric induced by $V(z)$, and $b^i(z)=\tilde{b}^i(z)-\frac{1}{2}\sum\limits_{l,j}(V^{-1}(z))^{lj}\Gamma_{lj}^i$ is the $i$ component of $b$. Here, $\Gamma_{lj}^i$ is the Christoffel symbols associated with this Riemannian metric, which satisfies 
\begin{equation}
\Gamma_{lj}^i=\frac{1}{2}\sum\limits_m g^{im}(\frac{\partial}{\partial{x^j}}g_{lm}+\frac{\partial}{\partial{x^l}}g_{jm}-\frac{\partial}{\partial{x^m}}g_{lj}),
\end{equation}
where $(g^{ij})(z)$ is the inverse of the Riemannian metric $(g_{ij})(z)=V(z)$. The divergence $\operatorname{div}b(z)$ is defined as
\begin{equation}
\mathtt{div}b(z)=\frac{1}{\sqrt{|V(z)|}}\sum\limits_i \frac{\partial}{\partial{z_i}}\left(b^{i}(z)\sqrt{|V(z)|}\right),
\end{equation}
where $|V(z)|$ is the determinate of Riemann metric. The Onsager-Machlup action functional can be considered as the integral of a Lagrangian 
\begin{equation}\label{lagrangian}
\begin{split}
L(z,\dot{z})=\frac{1}{2}[(\dot{z}-b(z))V(z)(\dot{z}-b(z))+\operatorname{div}b(z)-\frac{1}{6}R(z)].
\end{split}
\end{equation}

The minimum value of the action functional \eqref{action} means the largest probability of the solution paths on a small tube. We can investigate the most probable transition pathway by minimizing the Onsager-Machlup action functional. It captures sample paths of the largest probability around its neighbourhood.

To get the minimizer of the problem \eqref{action}, if one conduct by the variational principle, the most probable transition pathway between the two meta-stable states satisfies the classical Euler-Lagrange equation
\begin{equation}
\frac{d}{dt} \frac{\partial}{\partial {\dot{z}}}L(z,\dot{z})= \frac{\partial}{\partial{{z}}} L(z,\dot{z})
\end{equation}\label{EL_eqn}
with initial state $z(t_0)$ and final state $z(t_{f})$.  

A proof of the Euler-Lagrange equations, since it is not required for the rest of our discussions, can be found in the subject of calculus of variations \cite{gelfand2000calculus}. Chen et al \cite{chen2022data} have recently done a work of computing the most probable transition pathway for a carbon cycle system. Solving numerical solutions of Euler-Lagrange equations like Equation\eqref{EL_eqn} becomes tricky when faced with high-dimensional systems, although the shooting method is an efficient method for solving two-point boundary value problems. Optimal control is another way to look at calculus of variations problems, in that we view things in a dynamical nature. Concretely, we may re-parameterize the trajectories $x(t)$ in the form of a control form. 

Mathematically, the Lagrangian $L$ and the Hamiltonian $H$ are related via a construction known as the $Legendre$ $transformation$. Since this transformation is classical and finds applications in many diverse areas (optimization, geometry, physics), we now proceed to describe it. However, we will see that it does not quite provide the right point of view for our future developments, and it is included here mainly for historical reasons. 

\subsection{Problem Formulation and Assumptions for Detecting the Most Probable Transition Pathway}\label{Optimal_problem}
The study of optimal control theory actually originates from the classical theory of calculus of variations\cite{Liberzon2013Calculus}, starting from the seminal work of Euler and Lagrangian in the 1700s. These works culminated in the so-called $Lagrangian$ $Mechanics$ that reformulated $Newtonian$ $Mechanics$ in terms of extremal principles. In short, variational deals with optimization problems on ``curves", which can be described as an infinite dimensional extension of traditional optimization problems. 

\subsubsection{General Form of Optimal Control Problem}
And now, we will formulate precisely the optimal control problem in the general setting. We come to the situation where equality constrains are imposed on the admissible trajectories. Consider the most general ordinary differential equation

$$\dot{x}(t)=f(t, x(t), \theta(t)), \quad t \in\left[t_0, t_f\right], \quad x\left(t_0\right)=x_0$$

Where $t$ is the independent variable, $x(t)=(x_{1},...,x_{d})^{T} \in  \mathbb{R}^d$ is the (time-dependent) state variables, $\theta=(\theta_{1},...\theta_{m})^{T} \in \Theta \subset \mathbb{R}^m$ is the control, with $\Theta$ the control set. We will assume that the control set $\Theta$ is closed (but it need not be bounded). And $f$ is supposed to holds for the following conditions:
$\left\{\begin{array}{l}f(t, x, \theta) \text { is continuous in } t \text { and } \theta \text { for all } x \\ f(t, x, \theta) \text { is continuously differentiable in } x \text { for all } t, \theta\end{array}\right.$\\
The above two conditions are proposed to ensure that the equation() is well-posed. \cite{2003Introduction} 
\textbf{remark:} There are two crucial points: (i) $f$ is supposed to be differentiable with respect to $\theta$. (ii) $t \longmapsto \theta(t)$ is regular. Or to consider in general case $\theta$ to be a essentially bounded function of $t$.\\
\textbf{The Cost Functional}   Let us define the objective functional. We now consider functional of the following form
\begin{equation}\label{costfunc}
J[\theta]=\int_{t_0}^{t_f} L(t, x(t), \theta(t)) d t+\Phi\left(t_f, x\left(t_f\right)\right)
\end{equation}

Here $L: \mathbb{R} \times \mathbb{R}^d \times \Theta \rightarrow \mathbb{R}$ is called the running cost, and $\Phi: \mathbb{R} \times \mathbb{R}^d \rightarrow \mathbb{R}$ is called the terminal cost. Such type of optimal control problem we call it the $Bolza$ problem of optimal control.\\
\textbf{Remark:}  In $Mayer$ problem, we do not have running cost term, $(t_{f},x_{f})$ is the point where the total cost is minimized, and if this minimum is unique then we do not have to define the target set a prior.\cite{Liberzon2013Calculus}\\
\textbf{Admissible Target Set}\label{target set}  We denote the admissible target set by $S$, particularly our consideration is the transition phenomenon of a system between two metastable states case. The target sets $S=\left\{t_1\right\} \times\left\{x_1\right\}$ which corresponds to the most restricted case of a $fixed-time$ and $fixed-endpoint$ problem. Readers who are interested in this part can see\cite{Liberzon2013Calculus} for more types of target sets.\\

Now, we are ready to state the form of optimal control problem as follows:
\begin{equation}\label{controlproblem}
\left\{\begin{array}{l}
\inf _{\theta} J[\theta]=\int_{t_0}^{t_f} L(t, x(t), \theta(t)) d t+\Phi\left(t_f, x\left(t_f\right)\right) \\
\text { subject to } \\
\dot{x}(t)=f(t, x(t), \theta(t)), \quad t \in\left[t_0, t_f\right], \quad x\left(t_0\right)=x_0 .
\end{array}\right.
\end{equation}
Here $t_{0}$ is the initial time, $x_{0}$ is the initial state. According to optimal control theory, we often consider these two to be fixed. This is the case that will be considered in the rest of this paper. Let us motivate this approach in the context of an encouraging example.

 \subsection{From Variational to Optimal Conrtol}\label{convert}
We mainly consider a generalized stochastic differential equation in $\mathbb{R}^d$:
\begin{equation}
\begin{split}\label{2.10}
dX(t)=\tilde{b}(X(t))dt+\sigma(X(t))dB(t), t\in[0,t_{f}],
\end{split}
\end{equation}
with initial condition $X(0)=x_0 \in \mathbb{R}^d$, where $\tilde{b}:\mathbb{R}^d\rightarrow\mathbb{R}^d$ is a regular function, $\sigma:\mathbb{R}^d\rightarrow\mathbb{R}^{d\times k}$ is a $d\times k$ matrix-valued function, and $B$ is a Brownian motion in $\mathbb{R}^k$. The well-posedness of stochastic differential equation \eqref{2.10} has been widely investigated \cite{karatzas2012brownian}. Here we are interested in the transition phenomena between two metastable states (taking as the stable equilibrium states of the corresponding deterministic system) if there exist. 

Anyone who are interested in investigating transition phenomena, should estimate the probability of the solution pathways on a small tube. The Onsager-Machlup theory of stochastic dynamical system \eqref{generalSDE} gives an approximation of the probability \eqref{probability}. According to our analysis in Section \ref{OM}, from the point of view of variational calculus, the Onsager-Machlup action functional is shown as equation \eqref{action}. In order to detecting the most probable transition pathway in a dynamical system, we conclude this into a unconstrained minimization problem as follows:
\begin{equation}\label{minimize_problem}
    S_{t_{f}}^{BM}(z)=\inf\left\{\frac{1}{2}\int_{0}^{t_{f}}[(\dot{z}-b(z))V(z)(\dot{z}-b(z))+\operatorname{div}b(z)-\frac{1}{6}R(z)]dt \textrm{ for all } t \in [0,t_{f}]\right\}.
\end{equation}
Which can be simplified as (see Section\ref{OM} for a detailed explanation of the relevant parameters)
\begin{equation}
  S_{t_{f}}^{BM}(z)=\inf\left\{\int_{0}^{t_{f}} L(z,\dot{z}) dt\textrm{ for all } t \in [0,t_{f}]\right\},
\end{equation}
where $L(z,\dot{z})$ denotes the Lagrangian\eqref{lagrangian}.
Refer to the equation\eqref{main_problem} and equation\eqref{generalSDE}, the above problem can be reformulated as a deterministic optimal control problems with constraints to the state variable:
\begin{equation}\label{main_problem}
    \begin{cases}
    \underset{\theta \in \Theta }{\operatorname{minimize}} & \frac{1}{2}\int_{0}^{t_{f}} [\theta^{2}+\operatorname{div}b(z)-\frac{1}{6}R(z)]dt+\Phi(z(t_{f})), \\ 
    \text{ subject to } & \dot{z}(t)=\tilde{b}(z(t))+\sigma(z(t))\theta(t), \\
    &  {z}(0)=x_0.
    \end{cases}
\end{equation}
Here $z: [0,t_{f}] \to \mathbb{R}^d$ is the state and $\theta \in \Theta$ is the feedback control. Function $\Phi(x)$ is the terminal cost, and it is defined by 
\begin{equation}\label{terminal}
    \begin{cases}
    \Phi(x) = \frac{(x-x_{t_{f}})^2}{(x-x_{t_{f}})^2+1}. 
    \end{cases}
\end{equation}
It is clear that the variational approach \eqref{minimize_problem} would take us on a path of overly restrictive regularity assumptions. Instead, we would like to establish the Hamiltonian maximization property more directly, not via derivatives. As we will soon see, working with a richer perturbation family is crucial for obtaining sharper necessary conditions for optimality.

\bigskip
\subsection{Statement of Pontryagin’s Maximum Principle (PMP)}\label{PMP}
Pontryagin’s Maximum Principle is a hallmark result in optimal control theory and the calculus of variations. It supplys a necessary condition for optimality. The reason we considering this is that Pontryagin's Maximun Principle(PMP) generalizes the Euler Lagrange equations in highly nontrivial ways. Moreover, forms a nature bridge between optimal control theory and deep learning. A quite significant representative work can be refered in \cite{2017Maximum}.
We state the general form of Pontryagin's Maximum Principle. First we need some definitions.\\
\bigskip

\textbf{Hamiltonian}
$$
\begin{aligned}
&H: \mathbb{R} \times \mathbb{R}^d \times \mathbb{R}^d \times \Theta \rightarrow \mathbb{R}, \\
&H(t, x, p, \theta)=p^{\top} f(t, x, \theta)-L(t, x, \theta) .
\end{aligned}
$$

For a control $\boldsymbol{\theta}=\left\{\theta(t): t \in\left[t_0, t_f\right]\right\}$, we say it is admissable if $\theta(t) \in \Theta$ for all $t \in\left[t_0, t_f\right]$.  We can rewrite the cost functional \eqref{costfunc} in terms of the Hamiltonian: 
 \begin{equation}\label{costfunctional}
J(\theta)=\int_{t_0}^{t_f}(\langle p(t), \dot{x}(t)\rangle-H(t, x(t), \theta(t), p(t))) d t+\Phi\left(x\left(t_f\right)\right)
\end{equation}\\
Combined with equation\eqref{main_problem}, all we need to do is Hamiltonian maximization. Let $\theta^{*}(.)$ be an
optimal control, by which we presently mean that it provides a global minimum: $J(\theta^{*}) \le J(\theta)$
for all piecewise continuous controls $\theta$. Let us state the Pontryagin's Maximum Principle for effectively handling this. 

\textbf{Pontryagin's Maximum Principle}\label{pmp}\\
Let $\theta^{*}(t)$: $[t_{0},t_{f}] \mapsto \Theta $ be a bounded, measurable and admissable control that optimizes \eqref{main_problem} and $x^{*}(t)$: $[t_{0},t_{f}] \mapsto \mathbb {R}^{n}$ be its corresponding optimal state trajectory. $\Phi(t,(t_{f}))$ comes from equation \eqref{terminal} denotes the terminal cost. Then, there exsits a function $p^{*}$: $[t_{0},t_{f}] \mapsto  \mathbb {R}^{n}$ and a constant $p^{*}_{0} \le 0$ satisfying $(p^{*}_{0},p^{*}_{t}) \neq (0,0)$. And they have following properties:\\
1) The joint evolution of $x^{*}$ and $p^{*}$ are governed by the following Hamilton's canonical equation.
\begin{equation}\label{PMP_canonical}
\left\{\begin{array}{l}
\dot{x}^*(t)=\nabla_p H\left(t, x^*(t), p^*(t), \theta^*(t)\right), \quad x^*\left(t_0\right)=x_0 \\
\dot{p}^*(t)=-\nabla_x H\left(t, x^*(t), p^*(t), \theta^*(t)\right), \quad p^*\left(t_f\right)=-\nabla_x \Phi\left(x^*\left(t_f\right)\right) \\
\end{array}\right.
\end{equation}
2) For each fixed $t$, Hamiltonian has a global maximum at $\theta = \theta^{*}$, i.e.
\begin{equation}\label{PMP_H}
H\left(t, x^*(t), p^*(t), \theta^*(t)\right) \geq H\left(t, x^*(t), p^*(t), \theta^(t)\right) \\
\end{equation}
~~~~~~~~~~~~~~~~~~~~~~~~~~~~~~~~$\forall \theta \in \Theta$ and a.e.t $\in\left[t_0, t_f\right]$

Optimal control decisions must be taken everywhere alone the optimal trajectories. This is, we are devoted to find the global optimal solution of equation \eqref{PMP_canonical} and equation \eqref{PMP_H}. Moreover, we assume that the controls are not restricted by some practical considerations, or it will make the space of admissible trajectories cumbersome. Pontryagin's Maximum Principle consists of Equation \eqref{PMP_canonical}-\eqref{PMP_H}. Of course, there are many other forms of PMP that interested readers can refer to \cite{Liberzon2013Calculus} for more details. We will not step deeper here, because the control problem we will study in the rest of this paper is the fixed-time and fixed-endpoint case, as foreshadowing earlier. \\
\textbf{Remark:} \\
$\bullet$ The equation \eqref{PMP_canonical} is called the canonical equations. The equation(a) is $state$ $equation$ which describes the evolution of the state under the optimal control process. The equation(b) is called the $co-state$ $equation$ with $p^{*}$ beging the $co-state$. The role of co-state equation is to propogate back optimality condition and is the adjoint of variational equation.\\
$\bullet$ The maximizationcondition \ref{PMP_H} is the key point of the maximum principle. It requires that the optimal control $\theta^{*}$ must globally maximize the Hamiltonian. \\
\subsection{Computing the Most Probable Transition Pathway Via Pontryagin's Maximum Principle}
In this subsection, let us focus on the problem of form \eqref{main_problem} proposed at the end of subsection \ref{convert}. We first introduce variable $p(t)$, which denotes the co-state of the state $x(t)$ in the dynamical system \eqref{generalSDE}. With the help of \eqref{costfunctional}, we will reformulated problem \eqref{lagrangian} as:

\begin{equation}\label{costfunction}
\inf J(\theta)=\int_{0}^{t_f}(\langle p(t), \dot{x}(t)\rangle-H(t, x(t), \theta(t), p(t))) dt.
\end{equation}.\\
So the goal problem has now become a task of Hamiltonian maximization by the Pontryagin's Maximum Principle \ref{pmp}.

So far, we have successfully completed the transformation from variational method to optimal control. Instead of solving complex Euler-Langrange equation \eqref{EL_eqn}, we use the Pontryagin's maximum principle to solve optimization problem \eqref{main_problem}. Once we get the optimal feedback control $\theta^{*}$, we can compute the most probable transition pathway by the last two of equation \eqref{main_problem}.

\section{Numerical Algorithm for the Pontryagin’s Maximum Principle}\label{Algorithm}
Model-based deep learning for optimal control was introduced in\cite{J2015Deep} for stochastic systems. In this approach, the optimal control law is approximated using a neural network. The reformulated Pontryagin's maximum principle associated with problem \ref{generalSDE} is expressed as:
\begin{equation}\label{pmp_final}
\left\{\begin{array}{l}
\dot{x}^*(t)=f(t, x^{*}(t), \theta^{*}(t))=\tilde{b}(x^{*}(t))dt+\sigma(x^{*}(t))\theta^{*}, \quad x^*\left(t_0\right)=x_0 \\
\dot{p}^*(t)=-\nabla_x H\left(t, x^*(t), p^*(t), \theta^*(t)\right), \quad p^*\left(t_f\right)=-\nabla_x \Phi\left(x^*\left(t_f\right)\right) \\
H\left(t, x^*(t), p^*(t), \theta^*(t)\right) \geq H\left(t, x^*(t), p^*(t), \theta^(t)\right) \\
~\forall \theta \in \Theta  \text{and a.e.t} \in\left[t_0, t_f\right]
\end{array}\right.
\end{equation}

The above system is called an (extended) Hamiltonian system. It partially characterizes the optimality of the problem. In cases where certain convexity conditions are presented, the Hamiltonian system of the above form fully characterizes an optimal control. In what follows, if $(x^{*}, \theta^{*})$ is an optimal pair and $p^{*}$ is the corresponding adjoint function, then $(x^{*}, p^{*}, \theta^{*})$ 
will be called an optimal triple.

\subsection{Restatemenr of Pontryagin's Maximum Principle(PMP)}
The Euler-Lagrange equation is a necessary condition for optimality in the context of the Basic Calculus of Variations Problem, where the boundary points are fixed but the curves are otherwise unconstrained \cite{Liberzon2013Calculus}. Note that the boundary condition in equation \eqref{pmp_final} specifies the value of $p$ at the end of the interval $ [t_{0},t_{f}]$, i.e., it is a terminal condition rather than an initial condition. We know that the first-order condition of optimality \cite{Liberzon2013Calculus} implies $H_{\theta}\mid_{*}\equiv0$. To express it in a detailed form, that is
 \begin{equation}\label{first_order}
H_{\theta}\left(t, x^*(t), \theta^*(t), p^*(t)\right)=0 \quad \forall t \in\left[t_0, t_{f}\right] .
\end{equation}
The above equation implies that the Hamiltonian function $H(t,x^{*}(t),\theta^{*}(t),p^{*}(t))$ has a stationary point at $\theta^{*}(t)$ for all $t\in [t_{0},t_{f}]$. This is a Hamiltonian stationary condition. Readers can refer to \cite{Liberzon2013Calculus} to find the second-order condition of optimality is \textbf{Legendre-Clebsch condition}, and if this stationary point is an extremum, then it is a necessarily a \textbf{$maximum$}.

We emphasize that the PMP is only a necessary condition for optimality, hence there can be cases where solutions to the PMP is not actually globally optimal. Nevertheless, in practice the PMP is often strong enough to give good solution candidates, and when certain convexity assumptions are satisfied the PMP becomes sufficient \cite{1963Necessary}.

It seems that we have to deal with a two-point boundary value problem first. The two-point boundary value problems of the kind that we saw in the necessary conditions for optimality stated in PMP \eqref{PMP_canonical} can be solved numerically using the shooting or forward-backward sweep methods. Computational algorithms of solving such problem have been studied for decades with an extensive literature. Solvers exist in various programming languages and computing platforms such as MATLAB and Python.

\textbf{Remark:}

$\bullet$ One should always remember that the Pontryagin's Maximum Principle provides only a necessary condition for optimality. In other words, it help us to pick out optimal control candidates, each of which needs to be confirmed whether it is indeed the optimal.

$\bullet$ The reader should keep in mind that an optimal control may not be exist. But fortunately, most of the interesting problems does exist, because the conditions provided by Pontryagin's Maximum Principle is strong enough. The Hamiltonian maximization condition \eqref{PMP_H} is stronger than a typical first-order condition \eqref{first_order}.

$\bullet$ $H$ is not only stationary, but globally maximized at an optimal control, which is a much stronger statement if H is not concave. Moreover, the PMP makes minimal assumptions on the parameter space $\Theta$; the PMP holds even when $f$ is non-smooth with respect to $\theta$, \cite{2017Maximum} even when $\Theta$ is a discreet subset of $\mathbb {R}^{m}$.

$\bullet$ The maximum principle are in fact necessary for local optimality when closeness in the $(x,\theta)$-space is measured by the $0-form$ for $x$ and $\mathcal L_{1}\ norm$  for $\theta$. Even though, we stress that the Hamiltonian maximization condition remains global. A sufficient condition of optimality and corresponding proof are stated. Readers can refer to \cite[Theorm 2.5]{1999Stochastic} for details.

 \subsection{Backward Propagation Based Method of Successive Approximations} \label{algorithm}
 From this point onward, we derive and analyze algorithms entirely in continuous time, which allows us to characterize errors estimates and convergence in a more transparent fashion. Chernousko and Lyubushin \cite{F1982Method} put forward an iterative method based on alternating propagation and optimization steps: method of successive approximations (MSA). We will devise our algorithm based on this method. The main idea of continuous averaging algorithm is to average a series of auxiliary points in the iteration process, in which each iteration point is obtained by solving the auxiliary programming problem.
 
The Pontryagin's maximum principle consists of two major components: the forward-backward Hamiltonian dynamics and the maximization for the optimal parameters at each time. Note that the first equation of equation \eqref{pmp_final} is independent of co-state $p$, and the second is the most time-consuming step in deep learning. Compared with Frank Wolf algorithm, the advantage of this method is that in each iteration process, the iteration step obtained by solving the linear search problem is not needed, and the iteration step is predetermined, so the MSA algorithm is simple. However, due to the lack of adjustment in the iteration process, the convergence speed of this algorithm is slow in large road network. Since $\theta^{*}$ has no explicit solution is to be evaluated in numerical algorithm of maximization, the computation load is increase. For above reasons, we propose a successive approximation method based on back propagation, which is summarized in Algorithm 1.\\

\textbf{Algorithm:} Back Propagation Based Extended MSA~~~~~~~~~~~~~~~~~~~~~
  
   1  Give a initial guess for $\theta^{0} \in \Theta$;
   
   2  for \textbf{k}=0 to $iterations$, do\\
      solve $\dot{x}^*(t)=\nabla_p H\left(t, x^*(t), p^*(t), \theta^*(t)\right), \quad x^*\left(t_0\right)=x_0$ \\
      solve $\dot{p}^*(t)=-\nabla_x H\left(t, x^*(t), p^*(t), \theta^*(t)\right), \quad p^*\left(t_f\right)=-\nabla_x \Phi\left(x^*\left(t_f\right)\right) \\$
      train $u_{optimal}$ = $MLP(u)$;

   3  When the loss function converges, the operation stops. And $x^*$ is the most probable transition pathway we solved.\\
    
A challenge is that boundary value problem solvers are sensitive to initial guess. It tends to diverges with a bad initial guess $\theta^{0}$. The most urgent thing is to give modification to its divergent behavior. Li et al. \cite{2017Maximum} made a great contribution to this problem by doing the error estimation and convergence analysis for the Basic MSA.

\section{Applications of PMP for Computing the Most Probable Transition Pathway} \label{application}
\subsection{A Double-well System}
\begin{figure}[t]
    \centering   \includegraphics[width=0.5\linewidth]{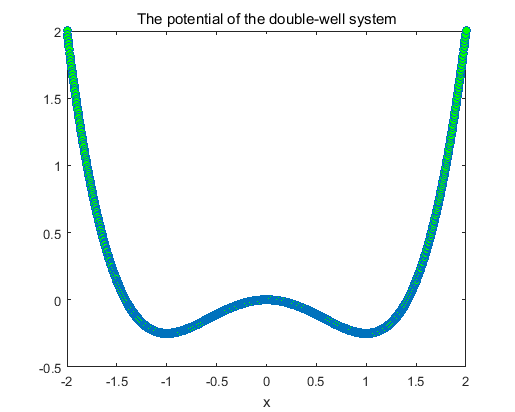}
    \caption{The potential of the system. And parameter $\sigma = 1$. There are two stable attractive basin: $(-1,0)$ and $(1,0)$. }\label{potential}
\end{figure}

\begin{figure}[htbp]
\begin{minipage}[]{0.5 \textwidth}
 \leftline{~~~~~~~\tiny\textbf{(\ref{fig:2}a)}}
\centerline{\includegraphics[width=6.5cm]{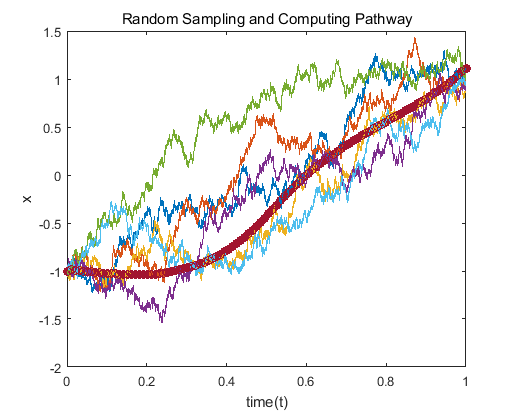}}
\end{minipage}
\hfill
\begin{minipage}[]{0.5 \textwidth}
 \leftline{~~~~~~~\tiny\textbf{(\ref{fig:2}b)}}
\centerline{\includegraphics[width=6.5cm]{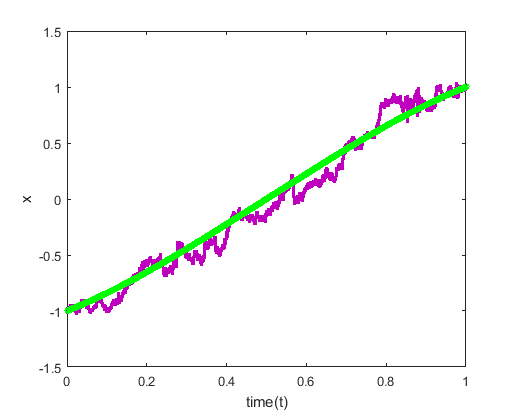}}
\end{minipage}
\caption{The most probable transition pathway and its verification. (2a):Random sampling of transition pathway and pathway calculated by Pontryagin's Maximum Principle. (2b): The most probable transition pathway calculated by Monte Carlo and Pontryagin's Maximum Principle, respectively. }\label{fig:2}
\end{figure}

To get your attention, our algorithm is validated with a simple one-dimensional model. Consider the following double-well system
\begin{equation}
\begin{split}
dx&=(x-x^3-)dt+ \sigma dW(t)\quad  t \in\left[0, T\right]\\
\end{split}
\end{equation}
Where $\sigma$ is the intensity of the noise, $W(t)$ is a wiener process. According to Fig. \ref{potential}, we consider the transitions between two metastable states under additive Brownian noise. Fig. \ref{fig:2} shows the results of numerical computation. It can be seen from Fig. \ref{fig:2}(a) that the transition pathway of numerically computing and that of random sampling are actually the same. Fig. \ref{fig:2}(a) well verifies the results of numerical computation, because we have made a comparison with the classical Monte Carlo simulation.

$$
$$
$$
$$

\subsection{Maier-Stein Sysytem}
Until now, we seem to have a certain intuitive sense of the numerical computation of Pontryagin's Maximum Principle. Overdamped systems without detailed balance arise in the study of disordered materials, chemical reactions far from equilibrium, and theoretical ecology.
As a typical example of the overdamped systems, the Maier-Stein system \ref{m_s} has been chosen to illustrate our results. We shall consider the following Maier-Stein system under multiplicative Brownian noise.

\begin{figure}[t]
    \centering
    \includegraphics[width=0.7\linewidth]{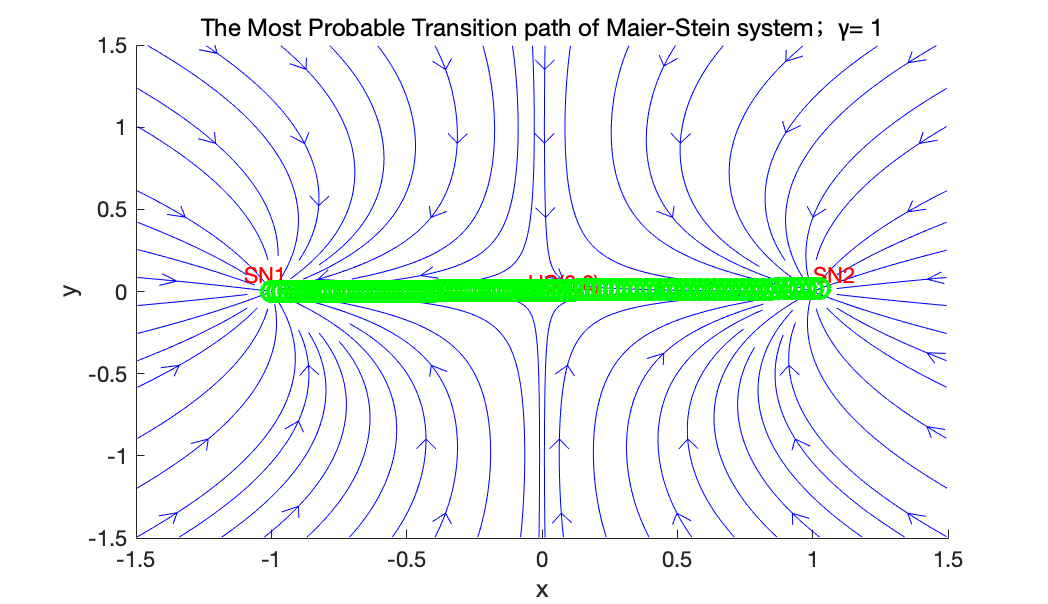}
    \caption{The most probable transition pathway of Maire-Stein system. The red dots indicate the equilibrium points, and the green pipes indicate the most probable transition pathway.}\label{gama_1}
\end{figure}

\begin{equation}\label{m_s}
\begin{split}
dx&=(x-x^3-\gamma xy^2)dt+\sigma_{1}(x, y) dW_{1}(t)\\
dy&=-(1+x^2)ydt+\sigma_{2}(x, y)dW_{2}(t)
\end{split}
\end{equation}

Where $W_{1}(t)$ and $W_{2}(t)$ are two independent noise variables, and we do not require them to follow the same distribution here. We study the most probable transition pathway between two metastable states of the system driven by multiplicative Brownian noise. In our parameter setting, $\sigma_{1}(x, y)$ =1+$\epsilon x^{2}$, where $\epsilon$ is a constant and $\sigma_{2}(x, y)$ = 1, which way we can get a diagonal diffusion matrix. According to the Pontryagin's Maximum Principle \ref{pmp_final}, and then we conclude as the following. And the Pontryagin's Maximum Princinple associated with Maire- Stein system is 
\begin{equation}\label{4.3}
\begin{aligned}
\dot{x}^*(t)&= x^*-(x^{*})^3-\gamma x^{*} (y^{*})^2 +(1+\epsilon (x^{*})^2) \theta_1^{*}\\
\dot{y}^*(t)&= -(1+(x^{*})^2) y^* + \theta_2^{*}\\
\dot{p_1}^*(t)&= p_1^{*} (3 (x^{*})^2+\gamma (y^{*})^2-2\epsilon x^* \theta_1^{*}-1) +2p_2 x^* y^* + \frac{\partial{L}}{\partial{x}}\\
\dot{p_2}^*(t)&= 2 p_1^{*} \gamma x^* y^* +p_2^{*} 
 (1+(x^*)^2)-2 \gamma y^*\\
\theta ^* &= argmax H\left(t, x^*(t), p^*(t), \theta^*(t)\right)\\
\end{aligned}
\end{equation} 

Solve the equation \eqref{4.3}, Fig. \ref{gama_1} visually show the most probable transition pathway in the sense of Onsager-Machlup functionals, which is consistent with the previous results obtained by the shooting method of \cite{2021Transition}. In fact, we can get more results in the sense of different $\gamma$, which avoids the influence of the existence of potential.

\subsection{Nutrient–Phytoplankton–Nooplanktonr System}
There are a lot of plankton in an aquatic (e.g. lake, river, sea) environment. Plankton are of two kinds, phytoplankton, or plant plankton, which photosynthesize and require essential elements (nutrient, e.g. nitrogen), and zooplankton, or animal plankton, which feed on phytoplankton. A remarkable feature associated with many phytoplankton population is the occurence of rapid and massive bloom formation \cite{2006Modelling}. When phytoplankton bloom occurs the number of phytoplankton sharply increases and decreases at short time and then returns its original low level.\\

It is important to explore the mechanism of phytoplankton bloom because of the important role of phytoplankton in the
earth. Plankton can affect large-scale global processes such as ocean-atmosphere dynamics and climate change. We introduce the following model to represent this evolutionary relationship.\\

\subsubsection{Model Introduction and Stability}
In this section, we will apply the Pontryagin's Maximum Principle method to the three-dimension dynamical system to compute the most probable transition pathway. For simplicity, the processes of recycling of phytoplankton and
zooplankton to nutrient are omitted just like that in paper \cite{2002A}. Let $N$ be the concentration of nutrient available for uptake, measured as mass per unit surface area of the ocean, $P$ the concentration of phytoplankton, and $Z$ the concentration of zooplankton. Then the model which describes the nutrient-phytoplankton-zooplankton dynamics is as follows
\begin{equation}\label{NPZ}
\left\{\begin{array}{l}
\frac{d N}{d t}=D\left(N_0-N\right)-f(N) P, \\
\frac{d P}{d t}=\alpha f(N) P-d_1 P-g(P) Z-w_1 P, \\
\frac{d Z}{d t}=\beta g(P) Z-d_2 Z-w_2 Z,
\end{array}\right.
\end{equation}
where
$$
\begin{aligned}
& f(N)= \begin{cases}\frac{b}{a} N, & 0 \leqslant N \leqslant a, \\
b, & N>a,\end{cases} \\
& g(P)=\frac{c P}{1+d P}.
\end{aligned}
$$
and all the parameters are positive. In equation \eqref{NPZ}, $N_{0}$ is constant input rate of nutrient, $D$, $w_{1}$ and $w_{1}$ are the washout rates for nutrient, phytoplankton and zooplankton respectively, $d_{1}$ and $d_{2}$ are the death rates for phytoplankton and zooplankton, $a$
and $b$ are the conversion rates from nutrient to phytoplankton and from phytoplankton to zooplankton. Let $d_1 + w_1 = D_1$, $d_2 + w_2 = D_2$. Then system \eqref{NPZ} becomes
\begin{equation}\label{NPZ_final}
\left\{\begin{array}{l}
\frac{d N}{d t}=D\left(N_0-N\right)-f(N) P, \\
\frac{d P}{d t}=\alpha f(N) P-g(P) Z-D_1 P \\
\frac{d Z}{d t}=\beta g(P) Z-D_2 Z .
\end{array}\right.
\end{equation}
The meaning of the parameter is described above.

\begin{figure}[t]
    \centering   \includegraphics[width=0.7\linewidth]{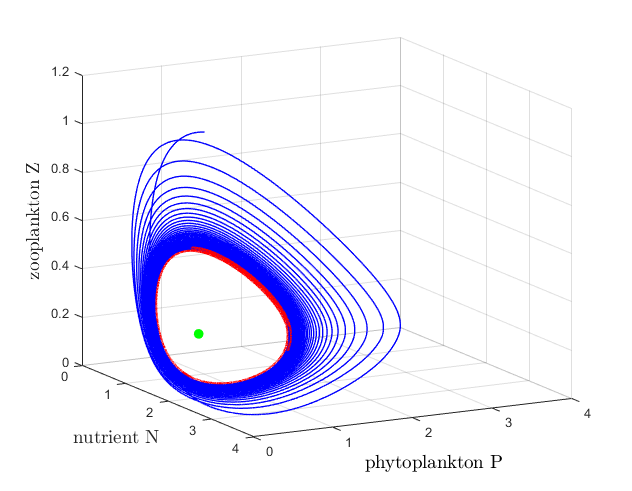}
    \caption{Bistable phenomenon: (1) a stable equilibrium (green point in the center) and a stable limit cycle (red curve). (2) The circular curve in blue is the asymptotic solution of this system. (3) The parameters are set as follows: $D = 0.1$, $a = 1$, $b = 1$, $\alpha= 1$, $c = 5$, $D_1 =0.2$, $\beta = 0.5$, $d = 0.1$, $D_2 = 2.1$, $N_0 = 9.96$.}\label{bistable}
\end{figure}

According to \cite{2012Hopf}, we set parameter values as the same, which results in bistable phenomenon in the system as in the Fig. \ref{bistable}. To further demonstrate the relationship between parameter $N_0$ and the dynamical behavior of the system, we leave all other parameters unchanged, as shown in Figure. However, Tianran Zhang did not take random fluctuations in an aquatic environment. There are a
lot of external disturbances in nature, such as climate change, and human activity. In order to explore the effect of noise disturbance on the dynamic
behavior of the Nutrient–Phytoplankton–Nooplan system, we consider a type of Gaussian noise purely from the perspective of mathematical model. We shall take into account the general random fluctuations rather
than infinitesimal random perturbations. To this end, the original system \eqref{NPZ_final} tends to be the following stochastic Nutrient–Phytoplankton–Nooplan model
\begin{equation}\label{SNPZ_final}
\left\{\begin{array}{l}
\frac{d N}{d t}=D\left(N_0-N\right)-f(N) P +\sigma_1 dW_t^1, \\
\frac{d P}{d t}=\alpha f(N) P-g(P) Z-D_1 P + \sigma_2 dW_t^2\\
\frac{d Z}{d t}=\beta g(P) Z-D_2 Z +\sigma_3 dW_t^3 .
\end{array}\right.
\end{equation}
Where $\sigma_i$ are the non-negative constant represents the intensity of the noise, $W_t^i$ are the Brownian motion.

\begin{figure}[t]
\begin{minipage}[]{0.5 \textwidth}

 \leftline{~~~~~\tiny\textbf{(\ref{fig:6}a)}}
\centerline{\includegraphics[width=6.5cm]{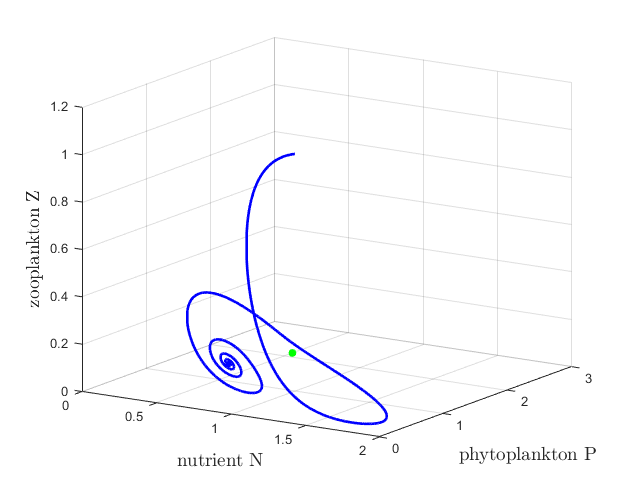}}
\end{minipage}
\hfill
\begin{minipage}[]{0.5 \textwidth}
 \leftline{~~~~~\tiny\textbf{(\ref{fig:6}b)}}
\centerline{\includegraphics[width=6.5cm]{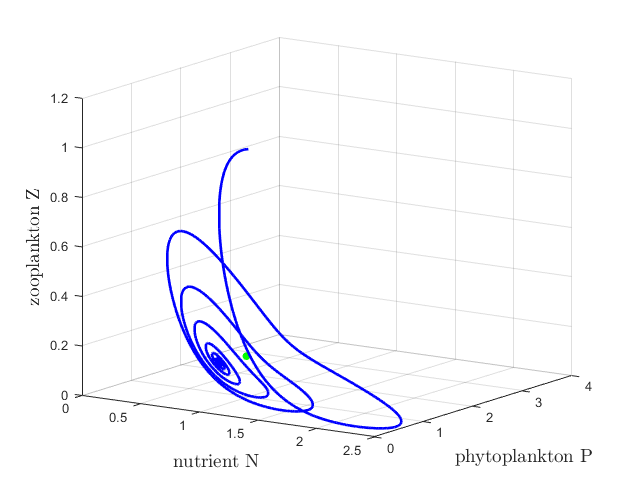}}
\end{minipage}
\begin{minipage}[]{0.5 \textwidth}
 \leftline{~~~~~\tiny\textbf{(\ref{fig:6}c)}}
\centerline{\includegraphics[width=6.5cm]{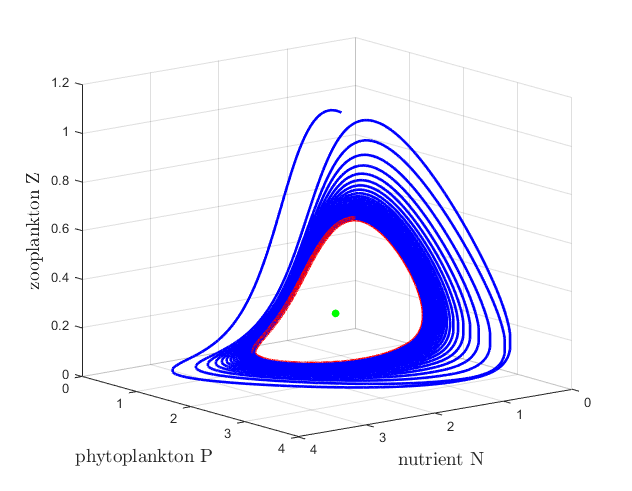}}
\end{minipage}
\hfill
\begin{minipage}[]{0.5 \textwidth}
 \leftline{~~~~~\tiny\textbf{(\ref{fig:6}d)}}
\centerline{\includegraphics[width=6.5cm]{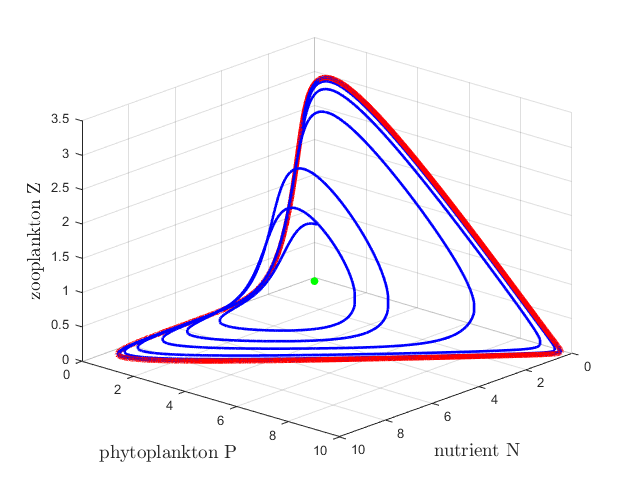}}
\end{minipage}
\caption{The dynamical diagram of different parameter value $N_0$. (6a) $N_0 = 6.0$; (6b) $N_0 = 8.0$; (6c) $N_0 = 10.0$; (6d) $N_0 = 15.0$. The green point is the coexisting equilibrium, the blue line is the solution pathway and the red circle is the stable limitcycle.}\label{fig:6}
\end{figure}

\subsubsection{Pontryagin's Maximum Principle for the N-P-Z System}
To better illustrate our approach, we first rearrange System \eqref{SNPZ_final} as the following problem
\begin{equation}\label{optimal_control_problem}
\left\{\begin{array}{l}
\inf _{\theta} J[\theta]=\int_{t_0}^{t_f} L(t, x(t), \theta(t)) d t+\Phi\left(t_f, x\left(t_f\right)\right) \\
\text { subject to } \\
\dot{x}(t)=f(t, x(t), \theta(t)), \quad t \in\left[t_0, t_f\right], \quad x\left(t_0\right)=x_0 .
\end{array}\right.
\end{equation}
where $L(t,x(t),\theta(t))$ stands for the Lagrangian of system \eqref{SNPZ_final}.  We take Onsager-Machlup functional theory according to additive Brownian noise, and 
\begin{equation}
\begin{split}\label{lagrange}
OM(X,\dot{X})=\frac{1}{2}\int_{0}^{t_{f}}[(\frac{\dot{X}-F}{\Sigma})^2+\operatorname{div}F(X)]dt,
\end{split}
\end{equation}
where $X=(N, P, Z)^T$, $\Sigma = (\sigma_1, \sigma_2, \sigma_3)^T$, the key transformation is that $\frac{\dot{X}-F}{\Sigma})^2 = W_{t}^{2}$, $W_t = (W_1, W_2, W_3)^2$. And now, we are ready to show the Pontryagin's Maximum Principle associated with problem \eqref{SNPZ_final} as follows

(i) $0 \leqslant x \leqslant a$
\begin{equation}
\begin{aligned}
\dot{x}^*(t)&= D\left(N_0-x^* \right)-\frac{b}{a}x^* y^* +\sigma_1 \theta_1^{*}\\
\dot{y}^*(t)&= \alpha \frac{b}{a}x^* y^*-\frac{c y^*}{1+d y^*} z^*-D_1 y^* +\sigma_2 \theta_2^{*}\\
\dot{z}^*(t)&= \beta \frac{c y^*}{1+d y^*} z^*-D_2 z^* +\sigma_3 \theta_3^{*} \\
\end{aligned}
\end{equation}

\begin{equation}
\begin{split}
\dot{p_1}^*(t)&= p_1^{*} (D+\frac{b}{a}y^*)+p_2^{*} \alpha\frac{b}{a}y^* +\frac{\alpha b}{2 a}  \\
\dot{p_2}^*(t)&= p_1^{*} \frac{b}{a}x^* + p_2^{*} (D_1 +z^* \frac{c}{(1+d y^*)^2} - \alpha \frac{b}{a}x^*) \\
&-p_3^{*} \beta z^* \frac{c}{(1+d y^*)^2} +\frac{1}{2}(\beta \frac{c}{(1+d y^*)^2} - \frac{b}{a} +\frac{2 z^* c d(1+d y^*)}{(1+d y^*)^4}) \\
\dot{p_3}^*(t)&= p_2^{*} \frac{c y^*}{1+d y^*} - p_3^{*} (\beta \frac{c y^*}{1+d y^*}-D_2) -\frac{c}{2 (1+d y^*)^2} \\
\end{split}
\end{equation}

(ii) $x > a$
\begin{equation}
\begin{aligned}
\dot{x}^*(t)&= D\left(N_0-x^* \right)-b y^* ++\sigma_1 \theta_1^{*}\\
\dot{y}^*(t)&= \alpha b y^*-\frac{c y^*}{1+d y^*} z^*-D_1 y^* +\sigma_2 \theta_2^{*}\\
\dot{z}^*(t)&= \beta \frac{c y^*}{1+d y^*} z^*-D_2 z^* +\sigma_3 \theta_3^{*}.\\
\end{aligned}
\end{equation}

\begin{equation}
\begin{aligned}
\dot{p_1}^*(t)&= p_1^{*} D \\
\dot{p_2}^*(t)&= p_1^{*} b + p_2^{*} (D_1 +z^* \frac{c}{(1+d y^*)^2} - \alpha b)\\
&-p_3^{*} \beta z^* \frac{c}{(1+d y^*)^2} +\frac{1}{2}(\beta \frac{c}{(1+d y^*)^2} +\frac{2 z^* c d(1+d y^*)}{(1+d y^*)^4}) \\
\dot{p_3}^*(t)&= p_2^{*} \frac{c y^{*}}{1+d y{*}} + p_3^{*} (D_2 - \beta \frac{c y^*}{1+d y^*}) -\frac{c}{2 (1+d y^*)^2} \\
\end{aligned}
\end{equation}

Please refer to the Appendix for more detaileds.

\subsubsection{Numerical Results and Analysis}
According to Section \ref{algorithm}, we use Python to compute the most probable transition pathway of the system in the sense of Onsager-Machlup functional.
\begin{figure}[t]
\begin{minipage}[]{0.5 \textwidth}
 \leftline{~~~~~\tiny\textbf{(\ref{fig:8}a)}}
\centerline{\includegraphics[width=6.5cm]{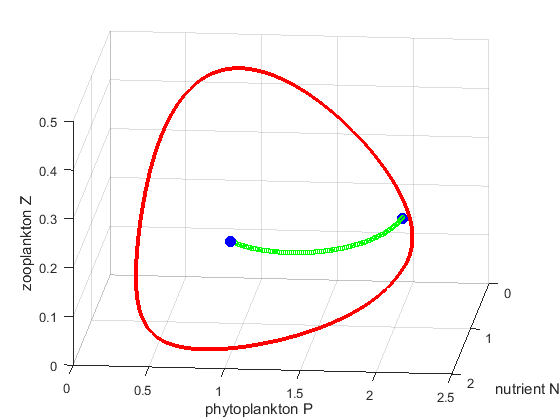}}
\end{minipage}
\hfill
\begin{minipage}[]{0.5 \textwidth}
 \leftline{~~~~~\tiny\textbf{(\ref{fig:8}b)}}
\centerline{\includegraphics[width=6.5cm]{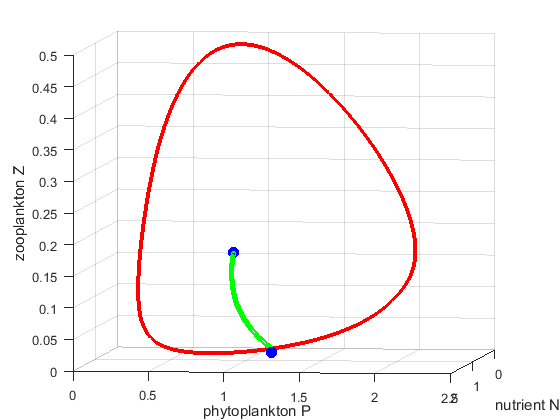}}
\end{minipage}
\begin{minipage}[]{0.5 \textwidth}
 \leftline{~~~~~\tiny\textbf{(\ref{fig:8}c)}}
\centerline{\includegraphics[width=6.5cm]{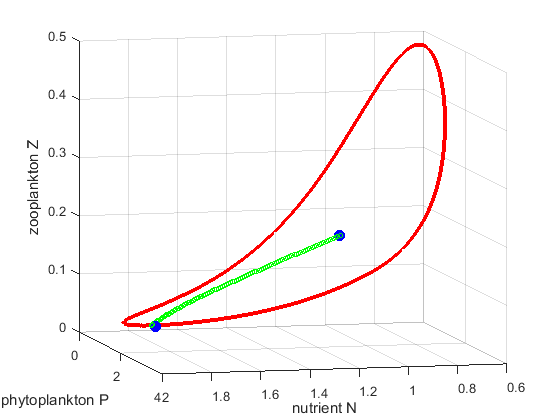}}
\end{minipage}
\hfill
\begin{minipage}[]{0.5 \textwidth}
 \leftline{~~~~~\tiny\textbf{(\ref{fig:8}d)}}
\centerline{\includegraphics[width=6.5cm]{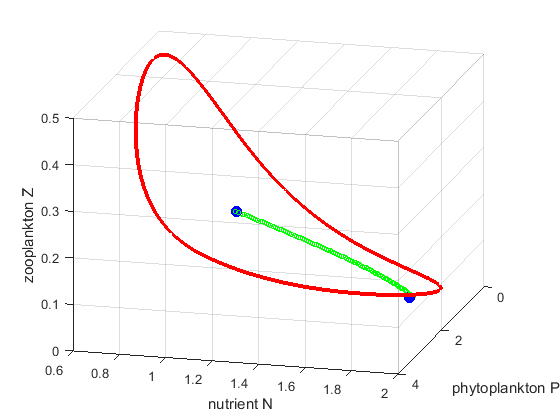}}
\end{minipage}
\caption{The most probable transition pathway in the sense of Onsager-Machlup. $T = 1$. The blue point is the coexisting equilibrium, the green line is the most probable transition pathway and the red circle is the stable limitcycle. }\label{fig:8}
\end{figure}

Fig.\ref{fig:8}(a) depicts the transition from the equilibrium point to the maximum $P$ value, where $P$ represents the phytoplankton concentration in the system. Fig. \ref{fig:8}(b - d) is the results of the maximum $N$ value, which is the state of maximum nutrient content. We provide three different angles for the most probable transition pathway. In the case of fixed time $T=1$, we know the influence of external random noise on the concentration of each variable in the system.
In this paper, the green transition pathway is not a real path, but represents the most probable tube. It is not difficult to find from the derivation of Onsager-Machlup functional that the green one is actually a tube with the highest probability, given a fixed $\epsilon$, it contains the largest number of sample transition pathways. \\ 

\textbf{Remark:} (a) We only provide limited results here, in fact we can solve for many other cases, and even extend the 3-dimensional case to higher dimensions, which will be explored in my future work. (b) We could also study the optimal transition time $T$, which contributes to the prevention of blooms in aquatic river systems. (c) Our calculation results can only be explained from the point of view of mathematical model, and may not have guiding significance for real aquatic river systems in nature, because we did not collect real experimental data. (d) Here we need to review the definition of reachable set( \ref{target set}), the transitions within the allowable error are accepted, there are a lot of inaccessible pathways because of the fixed transition time.

\section{Concludion and Discussion}\label{disscussion}
In fact, optimal control is widely used and has a mature theoretical basis. In this paper, we mainly study a kind of fixed-time and fixed-endpoint problem. There still remains other kinds of optimal control problem. In fact, the optimal control theory is an effective generalization of the variational. In the previous work, the most probable transition pathway of the dynamical system is studied from the perspective of functional analysis, combined with the OM functional principle, and obtained by solving the Euler-Lagrange equation. But this is tricky in higher dimensions. The optimal control theory shows great tolerance in the face of the diversity and complexity of noise perturbation. Considering the application of optimal control theory in dynamical system, we can consider solving more practical problems, such as the pursuit problem in ship ocean system, the hit problem of missile flight pair, and the problem of celestial body activity in planetary system.

The maximum principle established in this paper is not the most general.
One can also consider the possibility that the initial state $x_{0}$ is not fixed. We can impose separate constraints on $x_{0}$ and $x_{f}$ or, more generally, we can require they belong to some surface $\mathcal S_{2}$ in $\mathbb R^{2n}$. Maybe we can apply the maximum principle to manifolds.

\section*{Acknowledgements}
The work was done partially while the author was visiting the Center for Mathematical Sciences, Huazhong University of Science and Technology, Wuhan, China. This work was partly supported by NSFC grants 11771449.

\section*{Data Availability}
The data that support the findings of this study are openly available in GitHub.\\
\url{https://github.com/Cecilia-ChenJY/PMP-for-Optimal-Control}

\bibliographystyle{abbrv}
\bibliography{references}

\begin{appendices}
\section{Calculation of Pontryagin's Maximum Principle}
\subsection{Pontryagin's Maximun Principle of stochastic Maire-Stein System Driven by Multiplicative noise}
\end{appendices}\label{appA}
\begin{appendices}
For the stochastic  Maire-Stein system \eqref{m_s}, we shall derive the Pontryagin's Maximum Principle in detail. Denote the diffusion matrix $\sigma(x,y)=
\left[
\begin{matrix}
   \sigma_1(x,y)  &  0\\
   0  & \sigma_2(x,y)
\end{matrix}
\right]
$ and the vector field $\tilde{b}=(\tilde{b}^1,\tilde{b}^2)$, where 
\begin{align*}
\tilde{b}^1(x,y)&=x-x^3 -\gamma x y^2,\\
\tilde{b}^2(x,y)&=-(1+x^2) y.
\end{align*}
Thus, the Riemanian metric is $V(x,y)=(g_{ij})_{2\times2}(x,y)=
\left[
\begin{matrix}
   \frac{1}{(1+ \epsilon x^2)^2}  &  0\\
   0  & 1
\end{matrix}
\right]
$ and its inverse matrix is $(g^{ij})_{2\times2}(z)=
\left[
\begin{matrix}
   (1+ \epsilon x^2)^2  &  0\\
   0  & 1
\end{matrix}
\right]
$ . The determinate of $V$ is $|V(x,y)|=\frac{1}{ (1+\epsilon x^2) ^2         }$.

\noindent\textbf{Compute the Christoffel symbols}
\begin{align*}
 \Gamma_{11}^1&=\frac{1}{2}g^{11}\left(\frac{\partial}{\partial{x}}g_{11}+\frac{\partial}{\partial{x}}g_{11}-\frac{\partial}{\partial{x}}g_{11}\right) + \frac{1}{2}g^{12}\left(\frac{\partial}{\partial{x}}g_{12}+\frac{\partial}{\partial{x}}g_{21}-\frac{\partial}{\partial{x}}g_{11}\right)\\
  &=\frac{1}{2} (1+ \epsilon x^2)^2 (\frac {-4 \epsilon x}{(1+\epsilon x^2)^3}) +0\\
  &=-\frac{2 \epsilon x}{1+\epsilon x^2}\\
 \Gamma_{12}^1&=\frac{1}{2}g^{11}\left(\frac{\partial}{\partial{x}}g_{21}+\frac{\partial}{\partial{y}}g_{11}-\frac{\partial}{\partial{x}}g_{12}\right)+\frac{1}{2}g^{12}\left(\frac{\partial}{\partial{x}}g_{22}+\frac{\partial}{\partial{y}}g_{21}-\frac{\partial}{\partial{y}}g_{12}\right)\\
  &=0=\Gamma_{21}^1,\\
  \Gamma_{22}^1&=\frac{1}{2}g^{11}\left(\frac{\partial}{\partial{y}}g_{21}+\frac{\partial}{\partial{y}}g_{12}- \frac{\partial}{\partial{x}}g_{22}\right)+\frac{1}{2}g^{12}\left(\frac{\partial}{\partial{y}}g_{22}+\frac{\partial}{\partial{y}}g_{22}- \frac{\partial}{\partial{y}}g_{22}\right)\\
  &=0,\\
  \Gamma_{11}^2&=\frac{1}{2}g^{21}\left(\frac{\partial}{\partial{x}}g_{11}+\frac{\partial}{\partial{x}}g_{11}- \frac{\partial}{\partial{x}}g_{11}\right)+\frac{1}{2}g^{22}\left(\frac{\partial}{\partial{x}}g_{12}+\frac{\partial}{\partial{x}}g_{21}- \frac{\partial}{\partial{y}}g_{11}\right)\\
  &=-\frac{1}{2} \varepsilon^2\mu^2\cdot\partial{y}\left(\frac{1}{\varepsilon^2\mu^2\cdot f^2(x)}\right)\\
  &=0,\\
  \Gamma_{12}^2&=\frac{1}{2}g^{21}\left(\frac{\partial}{\partial{x}}g_{21}+\frac{\partial}{\partial{y}}g_{11}- \frac{\partial}{\partial{x}}g_{12}\right)+\frac{1}{2}g^{22}\left(\frac{\partial}{\partial{x}}g_{22}+\frac{\partial}{\partial{y}}g_{21}- \frac{\partial}{\partial{y}}g_{12}\right)\\
  &=0=\Gamma_{21}^2,\\
  \Gamma_{22}^2&=\frac{1}{2}g^{21}\left(\frac{\partial}{\partial{x}}g_{21}+\frac{\partial}{\partial{y}}g_{11}- \frac{\partial}{\partial{x}}g_{12}\right)+\frac{1}{2}g^{22}\left(\frac{\partial}{\partial{y}}g_{22}+\frac{\partial}{\partial{y}}g_{22}- \frac{\partial}{\partial{y}}g_{22}\right)\\
  &=0.
\end{align*}

\noindent\textbf{Compute the divergence}
The modified vector field is given by 
\begin{align*}
    b^1&=\tilde{b}^1-\frac{1}{2}\sum_{lj}(V^{-1}(z))^{lj} \Gamma^1_{lj}=\tilde{b}^1-\frac{1}{2} (1+\epsilon x^2)^2 \frac{-2 \epsilon x}{1+ \epsilon x^2}=\tilde{b}^1+\epsilon x (1+\epsilon x^2)\\
    &= x-x^3-\gamma x y^2+\epsilon x (1+ \epsilon x^2),\\
    b^2&=\tilde{b}^2-\frac{1}{2}\sum_{lj}(V^{-1}(z))^{lj} \Gamma^2_{lj}=\tilde{b}^2\\
    &= -(1+x^2) y.
\end{align*}
Therefore, the divergence of the modified vector field is
\begin{align*}
   \mathtt{div}b(x, y)&= (1+\epsilon x^2)(\frac {\partial}{\partial{x}} (b^1 \frac{1}{1+\epsilon x^2}) + \frac {\partial}{\partial{y}} (b^2 \frac{1}{1+\epsilon x^2}))\\
   &= -4x^2 - \gamma y^2 +\epsilon +3 \epsilon^2 x^2\\ 
   &+ \frac{-2\epsilon x(x-x^3-\gamma x y^2 +\epsilon x +\epsilon^2 x^3)}{1+\epsilon x^2}.\\
\end{align*}
\noindent\textbf{Compute the scalar curvature}
We use the notation of Einstein sum in the following computing. The scalar curvature is defined as the trace of the Ricci curvature tensor with respect to the Riemanian metric
\begin{align*}
R=g^{ij}R_{ij},
\end{align*}
where $R_{ij}$ is the Ricci curvature tensor. By \cite[Eqn 4.1.32, Eqn 4.3.6, Eqn 4.3.18]{jost2008riemannian}, we have
\begin{align*}
R&=g^{ik} R_{ik}=g^{ik}(g^{jl} R_{ijkl})=g^{ik}(g^{jl} (g_{im}R^{m}_{jkl}))\\
&=g^{ik}(g^{jl}(g_{im}(\frac{\partial}{\partial {x^k}}\Gamma_{lj}^m-\frac{\partial}{\partial{x^l}}\Gamma_{kj}^m+\Gamma_{ka}^m\Gamma_{lj}^a-\Gamma_{la}^m\Gamma_{kj}^a)).
\end{align*}
Since $g^{jl}g_{im}=\delta_{i=j,l=m}$, we have
\begin{align*}
R&=g^{ik}R^{m}_{ikm}=g^{ik}(\frac{\partial}{\partial{x^k}}\Gamma_{im}^m-\frac{\partial}{\partial{x^m}}\Gamma_{ki}^m+\Gamma_{ka}^m\Gamma_{mi}^a-\Gamma_{ma}^m\Gamma_{ki}^a).
\end{align*}

Since $\Gamma_{11}^1=-\frac{2 \epsilon x}{1+\epsilon x^2}$, $\Gamma_{jk}^i=0 (\text{others})$, $g^{11}= (1+\epsilon x^2)^2$, $g^{22}=1$, and $g^{12}=g^{21}=0$, we have
\begin{align*}
R=&g^{11}(\frac{\partial}{\partial{x}}\Gamma_{1m}^m- \frac{\partial}{\partial{x^m}}\Gamma_{11}^m+\Gamma_{1a}^m\Gamma_{m1}^a-\Gamma_{ma}^m\Gamma_{11}^a)\\
&+g^{22}(\frac{\partial}{\partial{y}}\Gamma_{2m}^m-\frac{\partial}{\partial{x^m}}\Gamma_{22}^m+\Gamma_{2a}^m\Gamma_{m2}^a-\Gamma_{ma}^m\Gamma_{22}^a)\\
=&g^{11}(\frac{\partial}{\partial{x}}\Gamma_{11}^1-\frac{\partial}{\partial{x}}\Gamma_{11}^1+\Gamma_{11}^1\Gamma_{11}^1-\Gamma_{11}^1\Gamma_{11}^1)\\
=&0.
\end{align*}

\noindent\textbf{Compute the Lagrangian}
By previous computing, the Lagrangian of the problem \eqref{minimize_problem} for the stochastic Maire-Stein system reduces to 
\begin{align*}
    L(z, \dot{z})=&(\dot{z}-b(z))V(z)(\dot{z}-b(z))+\mathtt{div}b(z)-\frac{1}{6}R(z)\\
                =& \Theta^2+\mathtt{div}b(z)-\frac{1}{6}R(z) \\
                =& \theta_{1}^{2}+ \theta_{2}^{2} +\mathtt{div}b(x, y) +\epsilon^2 x^2 -2 \epsilon x \theta_{1}.\\
\end{align*}
Where $z=(x,y)$, and $\mathtt{div}b(x, y)$ is calculated as above. Further we can get the expression for the Hamiltonian as
$$
H = p_1(x-x^3-\gamma x y^2 +(1 +\epsilon x^2) \theta_1)+ p_2( -(1+ x^2) y +\theta_2)-L.
$$

\noindent\textbf{Compute the Pontryagin's Maximum Principle}
\begin{enumerate}
    \item  
\begin{align*}
\dot{x}^*(t)&= x^*-(x^{*})^3-\gamma x^{*} (y^{*})^2 +(1+\epsilon (x^{*})^2) \theta_1^{*}\\
\dot{y}^*(t)&= -(1+(x^{*})^2) y^* + \theta_2^{*}\\
\end{align*}
\item
\begin{align*}
\dot{p_1}^*(t)&= p_1^{*} (3 (x^{*})^2+\gamma (y^{*})^2-2\epsilon x^* \theta_1^{*}-1) +2p_2 x^* y^* + \frac{\partial{L}}{\partial{x}}\\
\dot{p_2}^*(t)&= 2 p_1^{*} \gamma x^* y^* +p_2^{*} 
 (1+(x^*)^2)-2 \gamma y^*\\
\end{align*}
\end{enumerate}
\subsection{Pontryagin's Maximun Principle of NPZ System Driven by Additive Noise}
First of all, we write the stochastic three-dimensional model \eqref{NPZ_final} in a unified form as
\begin{align*}
\dot X& = F dt +\Sigma dW_t,\\
X& =(x, y, z)^T, F=(f_1, f_2, f_3)^T, W_t = (w_1, w_2, w_3)^T, \Sigma=(\sigma_1, \sigma_2, \sigma_3)^T.\\
\end{align*}
According to the theorem of Onsager-Machlup functional, we get the Lagrangian of system \eqref{NPZ_final} that is driven by additive noise as
\begin{align*}
L(X, \dot X)&=\frac{1}{2}((\frac{\dot X-F}{\Sigma})^2 +\mathtt{div}F(X)).\\
\mathtt{div}F(X)&= \frac{\partial{f_1}}{\partial{x}}+ \frac{\partial{f_2}}{\partial{y}} +\frac{\partial{f_3}}{\partial{z}}\\
&= -D - \frac{b}{a} y +\alpha \frac{b}{a} x -D_1 - z \frac{c}{(1+ d y)^2}+\beta \frac{c y }{1+d y} - D_2 \\
\end{align*}
for (i) $0 \leqslant x \leqslant a$. 
\begin{align*}
\mathtt{div}F(X)&= \frac{\partial{f_1}}{\partial{x}}+ \frac{\partial{f_2}}{\partial{y}} +\frac{\partial{f_3}}{\partial{z}}\\
&= -D +\alpha b -D_1 - z \frac{c}{(1+ d y)^2}+\beta \frac{c y }{1+d y} - D_2 \\
\end{align*}
for (ii) $x > a$.

And the Hamiltonian is $H = P^T \dot X -L(X, \dot X)$, where $P = (p_1, p_2, p_3)^T$. Then we get the Pontryagin's Maximum Principle.
\end{appendices}
\end{document}